\documentclass[11pt,reqno]{amsart}
\usepackage[T1]{fontenc}
\usepackage[english]{babel}
\usepackage[utf8]{inputenc}
\usepackage[square,numbers]{natbib}
\usepackage{amssymb}
\usepackage{amsmath}
\usepackage{xcolor}
\usepackage{amsfonts}
\usepackage{amsthm}
\usepackage{csquotes}
\usepackage{pgfplots}
\usepackage{tikz}
\usetikzlibrary{shapes,backgrounds}
\usetikzlibrary{arrows}
\usetikzlibrary{matrix}
\usetikzlibrary{cd}
\usepackage{latexsym}
\usepackage{enumerate}
\usepackage{graphicx}
\usepackage[a4paper,margin=1in]{geometry}
\usepackage{setspace}
\usepackage{array}
\usepackage{afterpage}
\usepackage{fontawesome5}
\usepackage{geometry}
\pgfplotsset{compat=1.17}

\usepackage{hyperref}
\setcitestyle{numbers}
\theoremstyle{definition}

\def\C{{\mathbb C}}
\def\R{{\mathbb R}}

\def\d{{\mathrm d}}
\def\C{{\mathbb C}}
\def\R{{\mathbb R}}

\def\C{{\mathbb C}}

\def\R{{\mathbb R}}

\def\d{{\mathrm d}}
\def\C{{\mathbb C}}
\def\R{{\mathbb R}}

\def\la{\lambda}

\begin{document}

\newtheorem{Theoreme}{Theorem}
\newtheorem{Corollaire}[Theoreme]{Corollary}
\newtheorem{Proposition}[Theoreme]{Proposition}
\newtheorem{Lemme}[Theoreme]{Lemma}
\newtheorem{Définition}[Theoreme]{Definition}
\newtheorem{Conjecture}[Theoreme]{Conjecture}
\newtheorem{Remarque}[Theoreme]{Remark}

 \title{ A note on Bézout type inequalities for mixed volumes and Minkowski sums.} 

\author{Cheikh Saliou Ndiaye  }
\thanks{This research is partly funded by the Bézout Labex, funded by ANR, reference ANR-10-LABX-58.}
\maketitle

\begin{abstract}
In this note, we  study Bézout type inequalities for mixed volumes and Minkowski sums of convex bodies in $\R^n$. We first give a new proof and we extend inequalities of  Jian Xiao on mixed discriminants. Then, we use the mass transport method to deduce some  Bézout type inequalities for mixed volumes. Finally, we apply these inequalities to obtain Bézout type inequalities for Minkowski sums.
 \end{abstract}
 \noindent{\bf Mathematics Subject Classification:} 52A39, 52A40.\\
{\bf Keywords:} Mixed volume, Minkowski sum, Convex body, optimal transport, Alexandrov-Fenchel, Bézout inequality.

\section{Introduction}

In 1779, Bézout gave the first version of his theorem about the number of intersection points of algebraic hypersurfaces: let $f_1,\dots,f_n$ be polynomials with coefficients in $\C^n$ such that the sets $X_i:=\{x\in \C^n | f_i(x)=0\}$ for $i\in \{1,\dots,n\}$ have no component in common, then 
\begin{equation}\#(X_1\cap\dots\cap X_n) \leq \deg(f_1)\dots\deg(f_n)\label{bz}.\end{equation}
 From the Bernstein-Kushnirenko-Khovanskii theorem (see \cite{Ber,Kho,Ku}), the quantities appearing in this inequality may be written in terms of mixed volumes of convex bodies. For any convex bodies $K_1,\dots,K_n$ in  $\mathbb{R}^n,$ their mixed volume is defined as 
  \begin{equation*}V(K_1,\dots,K_n)=\frac{1}{n!}\sum_{k=1}^n (-1)^{n+k} \sum_{i_1<\cdots<i_k\leq n} |K_{i_1}+\cdots+K_{i_k}|,\end{equation*}
  where $K_{i_1}+\cdots+K_{i_k}:=\{x_{i_1}+\cdots+x_{i_k}| x_{i_j} \in K_{i_j}\}$ is the Minkowski sum and $|\cdot|$ is the Lebesgue measure. Minkowski showed that for any convex bodies $K_1,\dots,K_m$ and non-negative $t_1,\dots,t_m \in \mathbb{R}_+$, the function $|t_1K_1+ \cdots +t_mK_m|$  can be extended as a polynomial function (see section 5.1 of \cite{SCH}):
\begin{equation*}\left|t_1K_1+ \cdots +t_mK_m\right|= \sum\limits_{|i|=n}\frac{n!}{i_1!i_2!\cdots i_m!}V(K_1[i_1],\dots,K_m[i_m])t_1^{i_1} \cdots t_m^{i_m},\end{equation*}
where $K[i_j]$ denotes  $K$ repeated $i_j$ times.
 Soprunov and Zvavitch  showed in \cite{SZ} that  the inequality \eqref{bz} can be rewritten in the following way: 
  for $ 2\leq r\leq n$, and for all choices of convex bodies  $P_1,\dots,P_r$, it holds 
\begin{equation} |\Delta|^{r-1}V(P_1,\dots,P_r, \Delta[n-r])\leq \prod_{i=1}^r V(P_i,\Delta[n-1]),\label{SopZva}\end{equation}
 where $\Delta$ is a $n-$dimensional simplex. Indeed, let $P_1,\dots,P_r$ be Newton polytopes and $f_1,\dots,f_r$ be polynomials having respectively $P_1,\dots,P_r$ as Newton polytopes and for $i=1,\dots,r, \\X_i:=\{x\in \C^n | f_i(x)=0\}$. Let $H_{r+1},\dots,H_n$ be linear forms, then  for all  $i\in \{r+1,\dots,n\}$,  the standard simplex conv$\{0,e_1,\dots,e_n\}$ is the Newton polytope of  $H_i$ and $\deg(H_i)=1$. The  Bernstein-Kushnirenko- Khovanskii theorem states that 
 $$\#(X_1\cap\dots\cap X_r\cap H_{r+1}\cap \dots\cap H_n)=n!V(P_1,\dots,P_r,\Delta[n-r])$$
 and  for all $i\in \{1,\dots,r\},\deg(H_i)=n!V(P_i,\Delta[n-1])$ (see \cite{JX,SZ} for more details). Therefore, in this case, the inequality \eqref{bz} becomes  \eqref{SopZva} by replacing $n!$ by $1/|\Delta|$. Soprunov and Zvavitch \cite{SZ} conjectured that if a convex body $A$ satisfies  $$|A|^{r-1}V(B_1,\dots,B_r, A[n-r])\leq \prod_{i=1}^r V(B_i,A[n-1]),$$ for all convex bodies  $B_1,\dots,B_r,$
 then $A$ is an $n-$dimensional simplex. Some positive partial answers about this conjecture are given in \cite{SSZ,SZ,MSz}. In \cite{SZ}, an  inequality  of the same flavor  where $A$ can be an arbitrary convex body was also studied and the question of  the best constant $b_{n,r}$ such that  for all convex bodies  $A,B_1,\dots,B_r$ in $\R^n$, we have
 \begin{equation}|A|^{r-1}V(B_1,\dots,B_r, A[n-r])\leq b_{n,r}\prod_{i=1}^r V(B_i,A[n-1])\label{tanich}\end{equation}
 was considered. Inspired by the works of Fradelizi, Giannopoulos, Hartzoulaki, Meyer, and Paouris in \cite{FGM,GHP}, it was proved in \cite{SZ} that $b_{n,r}\leq \frac{(nr)^r}{r!}$. Xiao used inequality \eqref{disc} below to show that $b_{n,r}\leq n^{r-1}$ in \cite{JX}. Brazitikos, Giannopoulos and Liakopoulos \cite{BGL} showed that $b_{n,r}\leq 2^{2^{r-1}-1}$. Our  first main result is the following improvement of the previously known bounds. We show that 
\begin{equation}b_{n,r}\leq\underset{k\in\{1,\dots,r\}}{\min}\left\{2^\frac{k(k-1)}{2}\frac{n^{r-k}}{(r-k)!}\right\}.\label{impXiao}\end{equation}
In particular, for $k=1$, this gives $b_{n,r}\leq \frac{n^{r-1}}{(r-1)!}$ and for $k=r$, we get $b_{n,r}\leq 2^\frac{r(r-1)}{2}$.
 Our main tools  improve Xiao's argument and the following Fenchel's inequality (\cite{Fen}, see also inequality (7.76) in \cite{SCH}):
\begin{multline}V(A[2],K_1,\dots,K_{n-2})V(B,C,K_1,\dots,K_{n-2})\\
\leq 2 V(A,B,K_1,\dots,K_{n-2}) V(A,C,K_1,\dots,K_{n-2}),\label{locAF} \end{multline}
for any convex bodies $A,B,C,K_1,\dots,K_{n-2}$ convex bodies in $\R^n$. We establish \eqref{impXiao} by using mixed discriminants. Let $M_1,\dots,M_m$ be positive semi-definite symmetric matrices in  $\mathbb{R}^n $ and $t_1,\dots,t_m \in \mathbb{R}_+$ then
\begin{equation*}
\det\left(t_1M_1+ \cdots +t_mM_m\right)= \sum\limits_{|i|=n}\frac{n!}{i_1!i_2!\cdots i_m!}D(M_1[i_1],\dots,M_m[i_m])t_1^{i_1} \cdots t_m^{i_m},\end{equation*}
 where the coefficient $D(M_1[i_1],\dots,M_m[i_m])$ is called the mixed discriminant. Mixed discriminants are non negative and for any  invertible map $T$, $D(TM_1,\dots,TM_n)$ is defined as
 \begin{equation} D(TM_1,\dots,TM_n):=|\det(T)|D(M_1,\dots,M_n),\label{factDet}\end{equation}
 (see section 5.5 of \cite{SCH}). In \cite{JX}, Xiao proved inequality \eqref{disc} which we extend  in theorem \ref{GX} below. Xiao's inequality states that for any integers $k,n$ such that $1\leq k\leq n$ and any  positive definite symmetric matrices $A,B,M_1,\dots,M_{n-k}$, we have 
\begin{equation}  \det(A)D(B[k],M_1,\cdots,M_{n-k})\leq \binom{n}{k}D(A[n-k],B[k])D(A[k],M_1,\cdots,K_{n-k}).\label{disc}\end{equation}

Our second main theorem is the application of this new Bézout type inequality for mixed volumes to establish a new Bézout type inequality for Minkowski sums. Let $c_{n,m}$ be the smallest constant such that for any convex bodies $A, B_1,\dots,B_m \subset \R^n$, 
\begin{equation*}|A|^{m-1}|A +B_1+\dots+B_m | \leq c_{n,m} \prod\limits_{k=1}^{m} |A+B_k|.\end{equation*}
Bobkov and Madiman  established in \cite{BM} that  $c_{n,m}\leq (m+1)^n$.  Fradelizi, Madiman  and Zvavitch  used a consequence of \eqref{disc} in \cite{FMZ} to show that  $\left( \frac{4}{3}+o(1)\right)^n \leq c_{n,2}\leq \varphi^n$ where $\varphi$ is the golden ratio. They also  gave more precise bounds for small dimensions: $c_{2,2}=1$  $c_{3,2}=4/3$ and they conjectured that  $c_{4,2}=3/2$. Several similar results for Minkowski sums of zonoids appear in \cite{FMMZ}.\\

In Section \ref{sec:Xiao}, we give an extension of inequality \eqref{disc} and then we use optimal transport as introduced in \cite{ADM} and developed in \cite{LX,JX} to deduce the same type  of inequality for mixed volumes. In Section \ref{bzt:MV}, we establish our Bézout inequality for mixed volumes \eqref{impXiao} and we conclude with some Bézout type inequalities for Minkowski sums in Section \ref{sec:Msum}.

  \noindent{\bf Acknowledgement:} We are grateful to Artem Zvavitch and Matthieu Fradelizi for introducing us to the question and their help. We also thank Pascal Dingoyan, Alfredo Hubard and Maud Szusterman for interesting discussions about this work.
  \section{Xiao type inequalities}\label{sec:Xiao}
Our first theorem is the following extension of inequality \eqref{disc} due to Xiao \cite{JX}.
\begin{Theoreme}\label{GX} Let $m\geq 1, n\geq 2 $ and $ i_1,\dots,i_m\geq 0$  be integers such that $|i|:=i_1+\dots+i_m\leq n$. Then, for any  positive semi-definite symmetric matrices $A,B_1,\dots,B_m,$ $ M_1,\dots,M_{n-|i|}$,
\begin{multline} \frac{n!}{i_1!\dots i_m!(n-|i|)!}\det(A)^{m}D(B_1[i_1],\dots,B_m[i_m],M_1,\dots,M_{n-|i|})\\
\leq \binom{n}{|i|}D(M_1,\dots,M_{n-|i|},A[|i|])\prod_{k=1}^m \binom{n}{i_k}D(B_k[i_k],A[n-i_k]).\label{VDXG}\end{multline}	
\end{Theoreme}
\noindent{\it Proof:} The proof is by induction on $m$. The case $m=1$ is  the inequality \eqref{disc} of  which we give a new proof inspired by ideas from \cite{AFO}.\\ 
\underline{Case $m=1$:} 
If $m=1,B_1=B$ and $i_1=k\leq n$, \eqref{VDXG} becomes  \eqref{disc}. Since $A$ and $B$ are positive semi-definite symmetric matrices, by simultaneous orthogonalization, there exists an invertible matrix $P$ and a diagonal matrix $\Lambda$ such that $A=PP^t$ and $B_1=P\Lambda P^t$. Using \eqref{factDet} it is enough to assume in \eqref{disc} that $A=I,B =\Lambda=$ diag$(\la_1,\dots,\la_n)$ and to prove that 
\begin{equation*}D\left(B[k],M_1,\dots,M_{n-k}\right)\leq \binom{n}{k}D\left(B[k],I[n-k]\right)D\left(I[k],M_1,\dots,M_{n-k}\right).\end{equation*}	
For any positive semi-definite symmetric matrices $C_1,\dots,C_n$, the  polarization formula of the mixed discriminant  is 
 \begin{equation*} D(C_1,\dots,C_n) =\frac{1}{n!}\sum_{\sigma \in S_n} \det\left(C_{\sigma(1)}^1,\dots,C_{\sigma(n)}^n\right),\end{equation*}
 where $C_i^j$ denotes the $j-th$ column of $C_i$ (see section 5.5 of \cite{SCH}). Taking $C_1=\dots=C_k=B$ and for $i \in \{k+1,\dots,n\},C_i=M_{i-k}$  and expanding the discriminants, one has
\begin{align*}
 D(B[k],M_1,\dots,M_{n-k})&= \frac{1}{n!}\sum_{|K|=k} k! \,\,\la_K\underset{\{\sigma(1),\dots,\sigma(k)\}= K}{\sum_{\sigma \in S_n}}  \det\left(\left(M_{\sigma(k+1)}^{K^c}\right)^1,\dots,\left(M_{\sigma(n)}^{K^c}\right)^{n-k}\right)\\
 &=\frac{1}{\binom{n}{k}}\sum_{|K|=k} \la_K \, D(M_{1}^{K^c},\dots,M_{n-k}^{K^c}),
\end{align*}
where for any $K \subset [n]:\{1,\dots,n\}$ such that $|K|=k$, $\la_K$ denotes $\prod_{j\in K}\la_j$, $K^c$ is the complementary of $K$ in $[n]$ and  $M^{K^c}$ is the $(n-k)$ dimensional matrix defined by $M^{K^c}=(M_{i,j})_{ i,j \in K^c}$ and then $D(M_{1}^{K^c},\dots,M_{n-k}^{K^c})$ is a mixed discriminant in dimension $n-k$. In particular, we have,
\begin{equation*} D(I[k],M_1,\dots,M_{n-k})=\frac{1}{\binom{n}{k}}\sum_{|K|=k}   D(M_{1}^{K^c},\dots,M_{n-k}^{K^c}) \end{equation*}
\begin{equation*} \text{and} \qquad  D(I[n-k],B[k])=\frac{1}{\binom{n}{k}}\sum_{|J|=k} \det(B^J)=\frac{1}{\binom{n}{k}}\sum_{|J|=k}  \la_J.
\end{equation*}
 It follows that
 \begin{align*}
\binom{n}{k}D(I[n-k],B[k])D(I[k],M_1,\dots,M_{n-k})&=\frac{1}{\binom{n}{k}}\sum_{|J|=k}  \la_J  \sum_{|K|=k}  D(M_1^{K^c},\dots,M_{n-k}^{K^c})\\
			&\geq \frac{1}{\binom{n}{k}}\sum_{|K|=k} \la_K D(M_1^{K^c},\dots,M_{n-k}^{K^c})\\
			&= D(B[k],M_1,\dots,M_{n-k}).
\end{align*} 
 Hence, the case $m=1$ is proved.\\
 \underline{Induction step: } Let us assume that inequality \eqref{VDXG} is true for $m-1$. We can  assume  that $A=I,B_m=$ diag$(\delta_1,\dots,\delta_n)$, then we have  
\begin{multline*}
 \frac{n!}{i_1!\dots i_{m}!(n-|i|)!} D\left(B_1[i_1],\dots,B_{m}[i_{m}],M_1,\dots,M_{n-|i|}\right)\\
 = \frac{n!}{i_1!\dots i_{m}!(n-|i|)!}\sum_{|J_m|=i_m}\frac{1}{\binom{n}{i_m}}\delta_{J_m}\quad D\left(B_1^{J_m^c}[i_1],\dots,B_{m-1}^{J_m^c}[i_{m-1}],M_1^{J_m^c},\dots,M_{n-|i|}^{J_m^c}\right)\\
 =\frac{(n-i_m)!}{i_1!\dots i_{m-1}!(n-|i|)!}\sum_{|J_m|=i_m}\delta_{J_m}\quad D\left(B_1^{J_m^c}[i_1],\dots,B_{m-1}^{J_m^c}[i_{m-1}],M_1^{J_m^c},\dots,M_{n-|i|}^{J_m^c}\right).
\end{multline*}
On the other hand, 
 \begin{multline*}
  \binom{n}{|i|}D(I[|i|],M_1,\dots,M_{n-|i|})\prod_{k=1}^m \binom{n}{i_k} D(I[n-i_k],B_k[i_k])\\
=\sum_{|J|=|i|}D(M_1^{J^c},\dots,M_{n-|i|}^{J^c})\prod_{k=1}^m \sum_{|J_k|=i_k} \det(B_k^{J_k})\\
  =\sum_{|J_m|=i_m}\delta_{J_m} \sum_{|J|=|i|}D(M_1^{J^c},\dots,M_{n-|i|}^{J^c}) \prod_{k=1}^{m-1} \sum_{|J_k|=i_k} \det(B_k^{J_k}).
\end{multline*}
Thus,  we need to prove that 
 \begin{multline*}
 \frac{(n-i_m)!}{i_1!\dots i_{m-1}!(n-|i|)!}\sum_{|J_m|=i_m}\delta_{J_m}\quad D\left(B_1^{J_m^c}[i_1],\dots,B_{m-1}^{J_m^c}[i_{m-1}],M_1^{J_m^c},\dots,M_{n-|i|}^{J_m^c}\right)\\
  \leq \sum_{|J_m|=i_m}\delta_{J_m} \sum_{|J|=|i|}D(M_1^{J^c},\dots,M_{n-|i|}^{J^c}) \prod_{k=1}^{m-1} \sum_{|J_k|=i_k} \det(B_k^{J_k}).
\end{multline*}

By comparing term by term, it is enough to prove that for each $J_m$,

 \begin{multline}\frac{(n-i_m)!}{i_1!\dots i_{m-1}!(n-|i|)!} D(B_1^{J_m^c}[i_1],\dots,B_{m-1}^{J_m^c}[i_{m-1}],M_1^{J_m^c},\dots,M_{n-|i|}^{J_m^c})\\
 \leq \sum_{|J|=|i|}D(M_1^{J^c},\dots,M_{n-|i|}^{J^c})\prod_{k=1}^{m-1}  \sum_{|J_k|=i_k} \det(B_k^{J_k})\label{sot}.
\end{multline}
 The induction hypothesis  allows to say that, for such a fixed $J_m$,
 \begin{multline*}\frac{(n-i_m)!}{i_1!\dots i_{m-1}!(n-|i|)!}D(B_1^{J_m^c}[i_1],\dots,B_{m-1}^{J_m^c}[i_{m-1}],M_1^{J_m^c},\dots,M_{n-|i|}^{J_m^c})\\
  \leq\underset{J^c \subset J_{m}^c}{\sum_{|J|=|i|}}D(M_1^{J^c},\dots,M_{n-|i|}^{J^c})\prod_{k=1}^{m-1} \underset{J_k \subset J_{m}^c}{ \sum_{|J_k|=i_k}} \det(B_k^{J_k}),
\end{multline*}
which implies \eqref{sot}. \hfill $\Box$
\begin{Remarque}

In \cite{LX} (Remark 3.6), the author affirmed that for any $|i|\leq n$,
$$\det(A)^{m-1}D\left(B_1[i_1],\dots,B_m[i_m],A[n-|i|]\right)\leq \frac{(n!)^{m-1}(n-|i|)!}{\prod_{j=1}^m (n-i_j)!}\prod_{k=1}^m D\left(B_k[i_k],A[n-i_k]\right).$$
is true if  $A,B_1,\dots,B_m$  are diagonal matrices and asked that whether it holds for  any positive semi-definite symmetric matrices. The induction method that we have used in the previous proof allows to answer positively to this question.
\end{Remarque}
\noindent\begin{Theoreme}\label{GXV}Let $m\geq 1, n\geq 2$  and $|i|:=i_1+\dots+i_m\leq n$ be integers. Then, for any   convex bodies $A,B_1,\dots,B_m$, $K_1,\dots,K_{n-|i|}$ in $\R^n$,
\begin{multline} \frac{n!}{i_1!\dots i_m!(n-|i|)!}|A|^{m}V(B_1[i_1],\dots,B_m[i_m],K_1,\dots,K_{n-|i|})\\
\leq \binom{n}{|i|}V(K_1,\dots,K_{n-|i|},A[|i|])\prod_{k=1}^m \binom{n}{i_k}V(B_k[i_k],A[n-i_k]).\label{bcm}\end{multline}	
\end{Theoreme}
 \noindent {\it Proof:} We follow the  method of Alesker, Dar and Milman \cite{ADM}, introduced in this context by Lehmann and Xiao \cite{LX}.
  Let $\gamma_n$ be the gaussian measure on $\R^n$ and for any convex body $K$, $\mathcal U(K)$ be  the uniform distribution on $K$. By  Brenier's theorem \cite{Br,Vil}, for $k=1,\dots,m,$ there exists a convex function $f_k$ such that $\nabla f_k$ pushes forward $\gamma_n$ onto $\mathcal U(B_k)$. Then $\nabla f_k(\R^n)=B$, 
 \begin{equation}|B_k|=\int_{\R^n} \det(\nabla^2 f_k),\label{volM}\end{equation}
and for $j=1,\dots,n-|i|$, there exists a convex function $h_j$ such that  $\nabla h_j$ pushes forward $\gamma_n$ onto $\mathcal U(K_j)$. It was proved \cite{ADM} that if $\nabla g_1$ and $\nabla g_2$ are two Brenier's maps, then 
\begin{equation*} \nabla g_1(\R^n)+\nabla g_2(\R^n)=(\nabla g_1+\nabla g_2)(\R^n).\end{equation*}
So, equality \eqref{volM}  also holds for Minkowski sums and therefore for mixed volumes. Hence
\begin{equation*}V(B_1[i_1],\dots,B_m[i_m],K_1,\dots,K_{n-|i|})=\int_{\R^n} \varphi(x)\d x,\end{equation*}
where $\varphi= D(\nabla^2f_1[i_1],\dots,\nabla^2f_{m}[i_m],\nabla^2 h_1,\dots,\nabla^2 h_{n-|i|})$. If $\int \varphi=0$ the theorem is proved, so, we assume that $\int \varphi\ne 0$. Let $\mu$ be the probability measure having density $\frac{\d\mu}{\d x}=\frac{1}{\int \varphi}\varphi(x)$. There exists a convex function  $f_A$ such that $\nabla f_A$ pushes forward $\mu$ onto $\mathcal U(A)$. Since Brenier's map satisfies the Monge-Ampère equation (see Subsection 4.1.1 of \cite{Vil}), one has
 \begin{equation}|A|\frac{\varphi(x)}{\int \varphi}= \det(\nabla^2 f_A(x))\quad \mu-\text{almost everywhere in } \R^n\label{MAmix}.\end{equation}
 Let $\psi=D(\nabla^2 f_A[|i|],\nabla^2 h_1,\dots,\nabla^2 h_{n-|i|})$, using Hölder inequality, we get that
\begin{multline*}
\binom{n}{|i|}V(K_1,\dots,K_{n-|i|},A[|i|])\prod_{k=1}^m \binom{n}{i_k}V(B_k[i_k],A[n-i_k])\\
= \binom{n}{|i|} \int_{\mathbb{R}^n}\psi(x)\d x\prod_{k=1}^m \int_{\mathbb{R}^n}\binom{n}{i_k}D(\nabla^2 f_k(x)[i_k] ,\nabla^2 f_A(x)[n-i_k])\d x\\
									\geq \left[\int_{\mathbb{R}^n}\left(\binom{n}{|i|}\psi(x)\prod_{k=1}^m \binom{n}{i_k}D(\nabla^2 f_k(x)[i_k] ,\nabla^2 f_A(x)[n-i_k])\right)^\frac{1}{m+1}\d x\right]^{m+1}.
\end{multline*}
The inequality \eqref{VDXG} in Theorem \ref{GX} allows to get 
\begin{multline*} \left[\int_{\mathbb{R}^n}\left(\binom{n}{|i|}\psi(x)\prod_{k=1}^m \binom{n}{i_k}D(\nabla^2 f_k(x)[i_k] ,\nabla^2 f_A(x)[n-i_k])\right)^\frac{1}{m+1}\d x\right]^{m+1}\\
\geq \left[\int_{\mathbb{R}^n}\left( \frac{n!}{i_1!\dots i_m!(n-|i|)!}\det\left(\nabla^2 f_A(x) \right)^{m} \varphi(x)\right)^{\frac{1}{m+1}}\d x\right]^{m+1}\\
 =\frac{n!}{i_1!\dots i_m!(n-|i|)!}V(B_1[i_1],\dots,B_m[i_m],K_1,\dots,K_{n-|i|}))|A|^{m}, \end{multline*}
by taking into account the equality \eqref{MAmix}.\hfill$\Box$\\

\begin{Corollaire}\label{BonItXiao}
Let $1\leq m\leq n$ , $|i|:=i_1+\dots+i_m\leq n$ be integers and $A,B_1,\dots,B_m$  be convex bodies in $\R^n$. Then, for any $j\in \{1,\dots,m \}$,
\begin{equation} \frac{n!(|i|-i_j)!}{i_1!\dots i_m!(n-i_j)!}|A|^{m-1}V(B_1[i_1],\dots,B_m[i_m], A[n-|i|])\leq\prod_{k=1}^m \binom{n}{i_k}V(B_k[i_k],A[n-i_k]).\label{vraiITxiao}\end{equation}	
\end{Corollaire}
\noindent{\it Proof:} Without loss of generality we can assume that $j=m$. Firstly we replace $m$ by $m-1$ in \eqref{bcm} (and maintain $|i|:=i_1+\dots+i_m$), then one has
\begin{multline}\frac{n!}{i_1!\dots i_{m-1}!(n-|i|+i_m )!}|A|^{m-1}V(B_1[i_1],\dots,B_{m-1}[i_{m-1}],K_1,\dots,K_{n-|i|+i_m})\\
\leq \binom{n}{|i|-i_m}V(K_1,\dots,K_{n-|i|+i_m},A[|i|-i_m])\prod_{k=1}^{m-1} \binom{n}{i_k}V(B_k[i_k],A[n-i_k]).\label{intt}\end{multline}
Taking $K_1=\dots=K_{i_m}=B_m$ in and $K_{i_m+1}=\dots=K_{n-|i|+i_m}=A$ in \eqref{intt} gives 
\begin{multline}\frac{(|i|-i_m)!}{i_1!\dots i_{m-1}!}|A|^{m-1}V(B_1[i_1],\dots,B_{m}[i_{m}],A[n-|i|])\\
\leq V(B_m[i_m],A[n-i_m])\prod_{k=1}^{m-1} \binom{n}{i_k}V(B_k[i_k],A[n-i_k]).\label{intt2}\end{multline}
We obtain the desired result by multiplying both sides of \eqref{intt2} by $\binom{n}{i_m}.$\hfill $\Box$
\begin{Remarque}
Notice that Corollary \ref{BonItXiao} improves  the following Xiao's inequality \cite{JX}: for any integers $2\leq m\leq n,$  $|i|:=i_1+\dots+i_m\leq n$  and convex bodies $A,B_1,\dots,B_m$ in $\R^n$, for any $j\in \{1,\dots,m \},$
\begin{equation*} \binom{n}{i_j}V(B_1[i_1],\dots,B_m[i_m], A[n-|i|])|A|^{m-1}\leq\prod_{k=1}^m \binom{n}{i_k}V(B_k[i_k],A[n-i_k]).\end{equation*}	
(Similar results are given in \cite{KH} with $A$ being the Euclidean ball.)
\end{Remarque}
\section{Bézout inequality for mixed volumes.}\label{bzt:MV}
In this section, we present upper bounds on the constant  $b_{n,r}$ defined in \eqref{tanich}. Here, our first argument is Lemma \ref{GFenc} below which is a generalization of Fenchel's inequality \eqref{locAF} and the second argument uses Corollary \ref{BonItXiao}.
 \begin{Theoreme}\label{luc}
Let $2\leq r\leq n$ be two integers  and $b_{n,r}$ be the best constant such that  for all convex bodies  $A,B_1,\dots,B_r$ in $\R^n$, it holds
 \begin{equation*}|A|^{r-1}V(B_1,\dots,B_r, A[n-r])\leq b_{n,r}\prod_{i=1}^r V(B_i,A[n-1]).\end{equation*}
 Then, $b_{n,r}\leq\underset{k\in\{1,\dots,r\}}{\min}\left\{2^\frac{k(k-1)}{2}\frac{n^{r-k}}{(r-k)!}\right\}.$
\end{Theoreme}
Before proving this theorem, we introduce the following generalization of Fenchel's inequality \eqref{locAF}.
\begin{Lemme}\label{GFenc}
Let $1\leq m\leq n$ be two integers  and  $A,B_1,\dots,B_m$ be convex bodies in $\R^n$. Then
 \begin{equation}|A|V(B_1,\dots,B_m, A[n-m])\leq 2^{m-1}V(B_1,\dots,B_{m-1}, A[n-m+1])V(B_m,A[n-1])\label{itAFloc}\end{equation}
 and 
 \begin{equation}|A|^{m-1}V(B_1,\dots,B_m, A[n-m])\leq 2^\frac{m(m-1)}{2}\prod_{i=1}^m V(B_i,A[n-1])\label{loclac}.\end{equation}
\end{Lemme}
\noindent{\it Proof:} We prove \eqref{itAFloc} by induction on $m$. The case $m=1$ is trivial and the case $m=2$ is Fenchel's inequality \eqref{locAF}. Let us  assume that \eqref{itAFloc} is verified for $m-1$. From Fenchel's inequality, 
\begin{multline}V(B_1,\dots,B_m, A[n-m])V(B_1,\dots,B_{m-2}, A[n-m+2])\\
\leq 2V(B_1,\dots,B_{m-1}, A[n-m+1])V(B_1,\dots,B_{m-2},B_m, A[n-m+1]).\label{endAF}\end{multline}
 By the induction hypothesis, we have
 \begin{multline}|A|V(B_1,\dots,B_{m-2},B_{m}, A[n-m+1])\\ \leq 2^{m-2}V(B_1,\dots,B_{m-2}, A[n-m+2])V(B_m,A[n-1]).\label{induc}\end{multline}
We end the proof by multiplying term by term \eqref{endAF} and \eqref{induc}.\\
 The proof of \eqref{loclac} follows from \eqref{itAFloc} by induction on $m$.
 \hfill $\Box$\\
\noindent{\it Proof of theorem \ref{luc}:} First we remark that, for any $k\in \{1,\dots,r\}$, by taking in \eqref{bcm}  $i_1=\dots=i_m=1,m=r-k$ and for $1\leq j\leq k, K_j=B_{r-k+j},$ for $k\leq j\leq n-|i|\},  K_{j}=A$, we get that 
 \begin{equation*}V(B_1,\dots,B_{r}, A[n-r])|A|^{r-k}\leq \frac{n^{r-k}}{(r-k)!}V(B_{r-k+1},\dots,B_{r}, A[n-k])\prod_{i=1}^{r-k} V(B_i,A[n-1]),\end{equation*}
According to \eqref{loclac}, we have
  \begin{equation*}V(B_{r-k+1},\dots,B_{r}, A[n-k])|A|^{k-1}\leq 2^\frac{k(k-1)}{2} \prod_{i=r-k+1}^{r} V(B_i,A[n-1]),\label{step_Xiao}\end{equation*}
  Hence the result follows. \hfill $\Box$

\section{Bézout inequality for Minkowski sums.}\label{sec:Msum}

In \cite{IZR}, Ruzsa showed that for any compacts sets $A,B_1,\dots,B_m \subseteq \R^n$ and any $\epsilon>0$, if $|A|\ne 0$, there exists a compact set $A'\subset A$ such that
\begin{equation}|A|^{m}|A'+B_1+\dots+B_m|\leq (1+\epsilon)|A'|  \prod_{k=1}^m|A+B_k|.\label{RUZSA}\end{equation}
It was noticed in \cite{FMZ} that \eqref{RUZSA} gives  
\begin{equation}|A|^{m-1}|B_1+\dots+B_m|\leq \prod_{k=1}^m |A+B_k|\label{rzsm}.\end{equation}
Notice that our methods allow to give a new simple proof of \eqref{rzsm}. Indeed, if $|i|=n$, inequality \eqref{bcm} becomes
\begin{equation} \frac{n!}{i_1!\dots i_m!}|A|^{m-1}V(B_1[i_1],\dots,B_m[i_m])\leq \prod_{k=1}^m \binom{n}{i_k}V(B_k[i_k],A[n-i_k])\label{mbot}\end{equation}
and since 

\begin{equation} |A|^{m-1}|B_1+\dots+B_m|=\sum_{i_1+\dots+i_m=n}\frac{n!}{i_1!\dots i_m!}|A|^{m-1}V(B_1[i_1],\dots,B_m[i_m]),\label{bingo}\end{equation}
\begin{equation} \prod_{k=1}^m |A+B_k|=\sum_{i_1=1}^n \cdots \sum_{i_m=1}^n \prod_{k=1}^m \binom{n}{i_k}V(B_k[i_k],A[n-i_k]).\label{binga}\end{equation}
According to \eqref{mbot}, each term of the right side of \eqref{bingo} is less than a corresponding term of the right side of \eqref{binga}. Hence, we get a new proof of \eqref{rzsm}. Another way of proving \eqref{rzsm} is to first show that for all $m\geq 1$ and any  positive semi-definite symmetric matrices $A,B_1,\dots,B_m$,  
\begin{equation*}\det(A)^{m-1}\det(B_1+\dots+B_m)\leq \prod_{k=1}^m \det(A+B_k), \end{equation*}
by extending each term into mixed discriminants like  in \eqref{bingo} and \eqref{binga} and applying Theorem \ref{GX}. Then we  use again the same optimal transport method.

\begin{Theoreme}\label{M:R} Let  $m, n\geq 1$, two integers and $c_{n,m}$ the best positive constant such that
\begin{equation}|A|^{m-1}|A +B_1+\dots+B_m | \leq c_{n,m}\prod\limits_{k=1}^{m} |A+B_k|,\label{intro3}\end{equation}
 for any convex bodies $A, B_1,\dots,B_m \subset \R^n$. Then 
 \begin{enumerate}[(i)]
 \item  \begin{equation*}  c_{n,m}\leq  \left(\frac{1-x_m}{1-mx_m}\right)^n< 2^n,\end{equation*}
 where $x_m$ is the unique real root of $P_r(x)=(1-mx)^m-(m-1)^{m-1}x^{m-1}(1-x)$.
\item\begin{equation*} c_{n,m} \geq \frac{e}{\sqrt{2\pi n(e-1)}}g(n,m)^{-(n+\frac{1}{2})},\end{equation*}
where
\begin{equation*} g(n,m)=\left[1-\frac{1}{m-1}\left(\frac{m-1}{m}\right)^{m}\right]^{m-1}.\end{equation*}
Furthermore, $\frac{4}{3}\leq g(n,m)< e^{e^{-1}}$
\end{enumerate}
\end{Theoreme}
\noindent Note that Theorem \ref{M:R} does not follow from \eqref{RUZSA} and the sets in \eqref{intro3} are convex bodies.
 \begin{Remarque}\label{RMKup} Here we give some  upper bounds of $c_{n,m}$ simpler but less optimal than that given in theorem \ref{M:R}
   \begin{enumerate}[(i)]
 \item
 for fixed $n$, $(c_{n,m})_{m\geq 2}$ is a increasing sequence with respect to $m$, 
$1 \leq c_{n,2} \leq c_{n,m} \leq c_{n,m+1} \leq  (c_{n,2})^m$ and  for $q>0$, we have   $c_{n,m+q} \leq c_{n,m}c_{n,q+1}$. In fact, in \eqref{intro3}, for $m\geq 3$, if $B_m=\{0\} $, we get $$  |A|^{m-2}\left|A +B_1+\dots+B_{m-1} \right| \leq \\c_{n,m}\prod\limits_{k=1}^{m-1} |A+B_k|,$$ hence, $c_{n,m-1} \leq c_{n,m}$. If $B_k= \{0\}$ for all $1\leq k\leq m$ then we obtain equality in \eqref{intro3} with $c_{n,m}=1$, so $c_{n,m}\geq 1$. It is shown in \cite{FMZ} that $1=c_{2,2} \leq c_{n,2}$ by using \eqref{locAF}. In \eqref{intro3}, if we replace $B_m$ by $B_m + \dots+B_{m+q}$, we get 
 \begin{equation}
 |A|^{m-1}\left|A + B_1+\dots+B_{m+q} \right| \leq c_{n,m}|A+B_m+\dots+B_{m+q}|\prod\limits_{k=1}^{m-1} |A+B_k|.\label{1hand}
\end{equation}
On the other hand, we have
 \begin{equation}
 |A|^{q}\left|A + B_m+\dots+B_{m+q} \right| \leq c_{n,q+1}\prod\limits_{k=m}^{m+q} |A+B_k|.\label{2hd}
\end{equation}
Thus, \eqref{1hand} and \eqref{2hd} allow to say that $c_{n,m+q}\leq c_{n,m}c_{n,q+1}$, in particular, $c_{n,m} \leq  (c_{n,2})^{m-1}$.\\
\item We also remark that a weaker bound for $c_{n,m}$ can be obtained by applying  \eqref{rzsm} to $\frac{1}{m}A + B_i$ instead of $B_i$ for $i \in \{1,\dots,m\}$ and $\left(1-\frac{1}{m}\right)A$ instead of $A$. This gives 
 \begin{equation*} |A|^{m-1}|A+B_1+\dots+B_m| \leq \left[\left(1+\frac{1}{m-1}\right)^{m-1}\right]^n\prod_{k=1}^m |A+B_k|,\end{equation*}
 so $c_{n,m}\leq \left[\left(1+\frac{1}{m-1}\right)^{m-1}\right]^n\leq e^n$.
 \item We extend again as sums both sides of \eqref{intro3} like \eqref{bingo} and \eqref{binga}  and compare each term of $|A|^{m-1}|A +B_1+\dots+B_m |$ by the term of the same index in $\prod\limits_{k=1}^{m} |A+B_k|$. It follows that $c_{n,m}\leq d_{n,m}$ where $d_{n,m}$ satisfies for all $|i|:=i_1+\dots+i_m\leq n$

 \begin{equation*} \frac{n!}{i_1!\dots i_m!(n-|i|)!}|A|^{m}V(A[n-|i|],B_1[i_1],\dots,B_m[i_m])\leq d_{n,m}\prod_{k=1}^m \binom{n}{i_k}V(B_k[i_k],A[n-i_k]).\end{equation*}
 In \eqref{bcm}, if $K_1=\dots=K_{n-|i|}=A$, it becomes
 \begin{equation*} \frac{n!}{i_1!\dots i_m!(n-|i|)!}|A|^{m-1}V(A[n-|i|],B_1[i_1],\dots,B_m[i_m])\leq \binom{n}{|i|}\prod_{k=1}^m \binom{n}{i_k}V(B_k[i_k],A[n-i_k]),\end{equation*}
 so any $d_{n,m}$ such that $ \binom{n}{|i|} \leq d_{n,m}$ for all $|i| \leq n $ gives an upper bound for $c_{n,m}$, hence 
 $$c_{n,m} \leq \underset{|i|\leq n}{\max} \binom{n}{|i|} = \binom{n}{\lfloor \frac{n}{2}\rfloor} < 2^n.$$
 \end{enumerate}
 \end{Remarque}
 \noindent{\it Proof of Theorem \ref{M:R} (i) (upper bound):} The following idea comes from \cite{FMZ} where the authors prove that $c_{n,2}\leq \varphi^n$ (where $\varphi$ is the golden ratio). Similarly to computations (iii) from Remark \ref{RMKup}, we recall that after expending   both sides of \eqref{intro3}   and comparing term by term, one has $c_{n,m}\leq d_{n,m}$ where $d_{n,m}$ satisfies for all $|i|:=i_1+\dots+i_m\leq n$
 \begin{equation*} \frac{n!}{i_1!\dots i_m!(n-|i|)!}|A|^{m}V(A[n-|i|],B_1[i_1],\dots,B_m[i_m])\leq d_{n,m}\prod_{k=1}^m \binom{n}{i_k}V(B_k[i_k],A[n-i_k]),\end{equation*}
thus
  \begin{equation*} |A|^{m}V(A[n-|i|],B_1[i_1],\dots,B_m[i_m])\leq \frac{i_1!\dots i_m!(n-|i|)!}{n!}d_{n,m}\prod_{k=1}^m \binom{n}{i_k}V(B_k[i_k],A[n-i_k]).\end{equation*}
  For any $j \in\{1,\dots,m\}$, \eqref{vraiITxiao} can be rewritten as 
  \begin{equation*} V(B_1[i_1],\dots,B_m[i_m], A[n-|i|])|A|^{m-1}\leq\frac{i_1!\dots i_m!(n-i_j)!} {(|i|-i_j)!n!}\prod_{k=1}^m \binom{n}{i_k}V(B_k[i_k],A[n-i_k]).\end{equation*}	
Therefore, for a fixed $j \in\{1,\dots,m\}$, to obtain an upper bound on $c_{n,m}$, it is enough to find $d_{n,m}$ such that for any $|i|\leq n$
$$\frac{i_1!\dots i_m!(n-i_j)!}{(|i|-i_j)!n!}\leq \frac{i_1!\dots i_m!(n-|i|)!}{n!}d_{n,m},$$
then, $\underset{|i|\leq n}{\max}\quad\binom{n-i_j}{|i|-i_j}\leq d_{n,m}$ and it follow that
$$c_{n,m} \leq \underset{|i|\leq n}{\max}\quad\underset{j\in\{1,\dots,m\}}{\min}\binom{n-i_j}{|i|-i_j}.$$
Note that $\underset{|i|\leq n}{\max}\quad\underset{j\in\{1,\dots,m\}}{\min}\binom{n-i_j}{|i|-i_j}\leq \underset{|i|\leq n}{\max}\binom{n}{|i|}<2^n$ as claimed in Remark \ref{RMKup}. For any $j \in\{1,\dots,m\}$, let $i_j=x_jn$ with $x_j\geq 0$ and $|x|:=x_1+\dots+x_m\leq 1$. According to Stirling formula, for any  integers $k\leq n$, one has $\binom{n}{k}\leq \frac{n^n}{(n-k)^{n-k}k^k}$, therefore:
$$c_{n,m} \leq \underset{|x|\leq 1}{\max}\quad\underset{j\in\{1,\dots,m\}}{\min}\left[\frac{(1-x_j)^{1-x_j}}{(1-|x|)^{1-|x|}(|x|-x_j)^{|x|-x_j}}\right]^n.$$
For $|x|\leq n$ and $ j\in\{1,\dots,m\}$, let
$$ F_j(x)=\frac{(1-x_j)^{1-x_j}}{(1-|x|)^{1-|x|}(|x|-x_j)^{|x|-x_j}}.$$
One has $\nabla^2 \log F_j(x)=\frac{-1}{1-|x|}UU^t$ where $U_j=\sqrt{\frac{|x|-x_j}{1-x_j}}$ and $U_k=1/U_j$ for $k\ne j$, then, $F_j$ is log-concave and $\min\log F_j$ also. Note that  $\underset{1\leq j\leq m}{\min} F_j$ is symmetric with respect to each hyperplane $\{x_i=x_k \}$, therefore $\underset{|x|\leq 1}{\max}\quad\underset{j\in\{1,\dots,m\}}{\min} F_j(x)$ is reached only if $x_1=\dots=x_m$. As a result, the problem becomes a simple study of a function having one variable:
 $$c_{n,m} \leq \underset{x\leq 1/m}{\max}\left[\frac{(1-x)^{1-x}}{(1-mx)^{1-mx}((m-1)x)^{(m-1)x}}\right]^n$$ and this maximum is reached at the unique $x_m \in]0,1/m[$ which is the real root of $P_r(x)=(1-mx)^m-(m-1)^{m-1}x^{m-1}(1-x)$. After some simplifications, we get $c_{n,m}\leq \left(\frac{1-x_m}{1-mx_m}\right)^n<2^n$. For $m=2$, $x_2=\frac{5-\sqrt{5}}{2}$ and $\frac{1-x_2}{1-2x_2}=\varphi$ (golden ratio), for $m=3$ we get that $c_{n,3}<(1.755)^n$.\hfill$\Box$\\

 For the lower bound of $c_{n,m}$, many results have been obtained in \cite{AAGJV,AFO,FGM,FMZ,FMMZ,SZ}. For the case $m=2$ we use methods from \cite{AAGJV,FMMZ,FMZ} and extend them to $m\geq 2 $. The idea is to represent volumes in \eqref{intro3} as volumes of a projection of $A$ to subspaces of lower  dimension and  optimize the lower bound with respect to those dimensions. We remark that \eqref{intro3} is equivalent to 
\begin{equation*}\frac{|A|^{m-1}|A + B_1+\dots+B_m |}{\prod\limits_{k=1}^{m} |A+B_k|} \leq c_{n,m}.\end{equation*}
So any choice of $A,B_1,\dots,B_m$ gives a lower bound of $c_{n,m}$. Let $\Gamma_2$ be the set of vector subspaces $\{E_1,E_2\}$ such that $\dim(E_1)=i,\dim(E_2)=j,0<i,j\leq n,i+j>n$. It is obtained  in \cite{AAGJV} that for any $A$ convex body in  $\R^n$,
\begin{equation} c_{n,2} \geq  \underset{A,\Gamma_2}{\max}\frac{|A|_n|P_{E_1\cap E_2} A|_{i+j-n}}{|P_{E_1} A|_i|P_{E_2} A|_j}=\frac{\binom{i}{i+j-n}\binom{j}{i+j-n}}{\binom{n}{i+j-n}}\label{DH}\end{equation}
where $A$ is a convex body in $\R^n$ and $P_E A$ is the projection of $A$ onto the subspace $E$. The equality case in \eqref{DH} is reached for a sharp choice of $E_1,E_2$ and $A$. The authors of \cite{FMZ} conclude that
\begin{equation*} c_{n,2}\geq \underset{i+j>n}{\underset{0<i,j\leq n }{\max} }\frac{\binom{i}{i+j-n}\binom{j}{i+j-n}}{\binom{n}{i+j-n}} \approx \frac{2}{\sqrt{\pi n}}\left(\frac{4}{3}\right)^n.\end{equation*}
By following the same idea, we get a lower bound for $c_{n,m}$.

 \begin{Theoreme}\label{hyp}
 Let $\alpha_1,\dots,\alpha_m \in \{1,\dots,n\}$ be integers and $d=|\alpha| -(m-1)n>0 $. Let $A\subset \R^n$ be a convex body and $\Gamma_m$ the set of collections of vector subspaces $\{E_1,\dots,E_m\}$ such that for $i=1,\dots,m$, $\dim(E_i)=\alpha_i$ and $E_i^\bot \subset \bigcap_{j \in [m] \smallsetminus \{i\}} E_j\}$. Then
  \begin{equation}c_{n,m} \geq \underset{\Gamma_m, A \text{ convex body }}{\max}
\frac{|A|_n|P_{\bigcap_{i =1}^m E_i} A|_d}{\prod_{i=1}^m |P_{E_i} A|_{\alpha_i}}.\label{proj}\end{equation}
\end{Theoreme}

\noindent{\it Proof:} Let $(e_1,\dots,e_n)$ be the canonical basis of $\R^n$. For $ i=1,\dots,m$, let $B_i= \sum_{k\in [n] \smallsetminus K_i}[0,e_k],$
 where $K_i= \{\alpha_1+\cdots+\alpha_{i-1}+ k |1\leq k \leq \alpha_i \}$. Thus, for any $i= 1,\dots,m$ and  any $t>0$, one has 
 $$|A+tB_i|_n \underset{t \rightarrow \infty}{\sim} t^{n-\alpha_i}|P_{E_i} A|_{\alpha_i}$$
  and
 \begin{equation*}
 |A+ tB_1+\dots+ tB_m|_n\underset{t \rightarrow \infty}{\sim} t^{nm-|\alpha|}\left|P_{\bigcap_{i=1}^m B_i} A\right|_{d}.
 \end{equation*}
Let $t$ to infinity,  then we get 
 \begin{equation*} |A|_n^{m-1}|P_{\bigcap_{i=1}^m E_i} A|_{d} \leq c_{n,m}\prod_{i=1}^m |P_{E_i}A|_{\alpha_i}.\end{equation*}
 The desired result follows. \hfill$\Box$\\

  \begin{Corollaire}\label{cor-alphaM}
 Let $\alpha_1,\dots,\alpha_m \in [1,n]$ be integers and $d=|\alpha| -(m-1)n>0 $. Then
  \begin{equation}
 c_{n,m} \geq \underset{\alpha}{\max}\frac{\prod_{i=1}^m \binom{\alpha_i}{d}}{\binom{n}{d}}.\label{bin}
\end{equation}
\end{Corollaire}
 \noindent{\it Proof:} Let $B_\infty(d)= \sum_{k=n-d+1}^n [-e_k,e_k]$. For $ i=1,\dots,m$, let $$J_i=[n-\alpha_1+\ldots+n-\alpha_{i-1}+1\quad;\quad n -\alpha_1+\ldots+n-\alpha_{i}],$$
$B_1(i)= \underset{J_i}{\text{conv}} (\pm e_k),$ and $E_i=\text{vect}\{e_k,k \in  [n]\smallsetminus J_i\}.$ If $A=\text{conv} \left\{\sum_{i=1}^m B_1(i),B_{\infty}(d) \right\}$, we get that  $$|A|=\prod_{i=1}^m |B_1(i) |\times|B_\infty(d)|/\binom{n}{d}\quad ;\quad |P_{\bigcap_{i=1}^m E_i}A |=|B_\infty(d)|$$ and for $ i=1,\dots,m$
$$|P_{E_i}A |= \prod_{k\in m\smallsetminus\{i\}}|B_1(k) |\times|B_\infty(d)|/\binom{\alpha_i}{d}.$$
Thus simplifying the volumes we get,
 \begin{equation*} \frac{ |A|_n^{m-1}|P_{\bigcap_{i=1}^m E_i} A|_{d} }{\prod_{i=1}^m |P_{E_i}A|_{\alpha_i}} =\frac{\prod_{i=1}^m \binom{\alpha_i}{d}}{\binom{n}{d}^{m-1}},\end{equation*}
 for all $\alpha$ satisfying the given hypothesis. Hence the result follows. \hfill $\Box\qquad$\\
\noindent {\it Proof of Theorem \ref{M:R} (ii) (lower bound):}  For the lower bound, Corollary \ref{cor-alphaM} allows to say that 
 $$c_{n,m}\geq \underset{d=|\alpha| -(m-1)n>0}{\underset{\alpha >0}{\max}} \frac{\prod_{i=1}^m \binom{\alpha_i}{d}}{\binom{n}{d}^{m-1}}.$$

 Now, using the Stirling formula, one has
 \begin{equation*}
 \frac{\prod_{i=1}^m \binom{\alpha_i}{d}}{\binom{n}{d}^{m-1}}\geq  \frac{1}{\sqrt{2\pi d}}\prod_{i=1}^m\sqrt{\frac{\alpha_i}{\alpha_i-d}}\times\frac{1}{d^d} \left[\sqrt{\frac{n-d}{n}}\frac{ (n-d)^{n-d}}{ n^n} \right]^{m-1}\prod_{i=1}^m\frac{ \alpha_i^{\alpha_i}}{ (\alpha_i-d)^{\alpha_i-d}}e^{\text o(1)}.
\end{equation*}
 For $i=1,\dots,m$, one assume that  $\alpha_i=nx_i$ and $d=ny$, then $y=|x|-(m-1)$ and 
\begin{equation*}\underset{\alpha}{\max}\left( \frac{1}{d^d} \left[\frac{ (n-d)^{n-d}}{ n^n} \right]^{m-1}\prod_{i=1}^m\frac{ \alpha_i^{\alpha_i}}{ (\alpha_i-d)^{\alpha_i-d}}\right)=\underset{x}{\max}f(x)^n,
\end{equation*}
where 
 \begin{equation*}f(x_1,\dots,x_m)= \frac{(1-y)^{(1-y)(m-1)} }{y^y}\prod_{i=1}^m\frac{ x_i^{x_i}}{ (x_i-y)^{x_i-y}}.
\end{equation*}

By studying this function, we easily find that $\max f=\frac{x^m}{y}>1$ where $$x=(m-1)\left[ m-\left(\frac{m-1}{m}\right)^{m-1} \right]^{-1}$$ and $y =m(x-1) +1$.  It follows that
\begin{align*}
 c_{n,m}&\geq\max \left(\frac{x^m}{y}\right)^n\frac{1}{\sqrt{2\pi ny}}(\sqrt{1-y})^{m-1}\left(\sqrt{\frac{x}{x-y}}\right)^m (1+\text{o}(1))\\
 		&=\frac{1}{\sqrt{2\pi n}} \left(\frac{x^m}{y}\right)^{n+1/2}\sqrt{\frac{(1-y)^{m-1}}{(x-y)^m}}(1+\text{o}(1)).\\
 \end{align*}
  
 Besides, $\frac{x^m}{y}=\left[1-\frac{1}{m-1}\left(\frac{m-1}{m}\right)^{m}\right]^{1-m}$ is a increasing sequence with respect to $m$ while \begin{equation*}\frac{(1-y)^{m-1}}{(x-y)^m}=\frac{\left(\frac{m}{m-1}\right)^{m-1}-\frac{1}{m}}{1-\frac{1}{m}-\left(\frac{m-1}{m}\right)^{m}}\quad\text{ is decreasing.}\end{equation*} 
 Thus
 \begin{equation*}\quad\frac{e}{\sqrt{e-1}}\left(\frac{4}{3}\right)^{n+1/2}< \left(\frac{x^m}{y}\right)^{n+1/2}\sqrt{\frac{(1-y)^{m-1}}{(x-y)^m}} < \sqrt6(e^{e^{-1}})^{n+1/2},\end{equation*}
 it leads to
 \begin{equation*} c_{n,m}\geq\frac{1}{\sqrt{2\pi n}}\frac{e}{\sqrt{e-1}}\left(\frac{4}{3}\right)^{n+1/2}.\end{equation*}
For $m=2$, $c_{n,2}\geq \frac{2}{\sqrt{\pi n}}\left(\frac{4}{3}\right)^n(1+\text{o}(1))$ which was found in \cite{FMZ}. For $m=3$, we find 
\begin{equation*}c(n,3)\geq \frac{81}{4\sqrt{115}\sqrt{\pi n}}\left(\frac{729}{529}\right)^n(1+\text{o}(1))\approx \frac{1.89}{\sqrt{\pi n}}(1.378)^n.\end{equation*}


  \faAt: \href{mailto:cheikh-saliou.ndiaye@univ-eiffel.fr}{cheikh-saliou.ndiaye@univ-eiffel.fr}. \faMapPin: Univ Gustave Eiffel, Univ Paris Est Creteil, CNRS, LAMA UMR8050 F-77447 Marne-la-Vallée, France
 \end{document}